\documentclass[12pt,a4paper,centertags]{amsart}
\usepackage{pslatex}
\usepackage{hyperref}
\newtheorem{theorem}{Theorem}
\newtheorem{lemma}[theorem]{Lemma}
\newtheorem{proposition}[theorem]{Proposition}
\newtheorem{corollary}[theorem]{Corollary}
\theoremstyle{definition}
\newtheorem{definition}[theorem]{Definition}
\theoremstyle{remark}
\newtheorem{remark}[theorem]{Remark}
\numberwithin{equation}{section}

\hypersetup{pdftitle={Regularity of transition semigroups associated to a 3D stochastic Navier-Stokes equation},
pdfsubject={A 3D stochastic Navier-Stokes equation with a suitable non degenerate additive
noise is considered. The regularity in the initial conditions of every Markov
transition kernel associated to the equation is studied by a simple direct
approach.},
pdfauthor={F. Flandoli, M. Romito},
pdfkeywords={stochastic Navier-Stokes equations, Markov selections, strong Feller property}}
\begin{document}
\title[Regularity of transition semigroups for the NSeqs]
    {Regularity of transition semigroups associated to a 3D stochastic Navier-Stokes equation}
\author[F. Flandoli]{Franco Flandoli}
\address{Dipartimento di Matematica Applicata, Universit\`a di Pisa,
         via Bonanno Pisano 25/b, 56126 Pisa, Italia}
\email{flandoli@dma.unipi.it}
\author[M. Romito]{Marco Romito}
\address{Dipartimento di Matematica, Universit\`a di Firenze,
         Viale Morgagni 67/a, 50134 Firenze, Italia}
\email{romito@math.unifi.it}
\date{}
\begin{abstract}
A 3D stochastic Navier-Stokes equation with a suitable non degenerate additive
noise is considered. The regularity in the initial conditions of every Markov
transition kernel associated to the equation is studied by a simple direct
approach. A by-product of the technique is the equivalence of all transition
probabilities associated to every Markov transition kernel.
\end{abstract}
\maketitle
\section{Introduction}
An old dream in stochastic fluid dynamics is to prove the well posedness of a
stochastic version of the 3D Navier-Stokes equations, taking advantage of the
noise, as one can do for finite dimensional stochastic equations with non
regular drift (see for instance Stroock \& Varadhan \cite{SV}). The problem is
still open, although some intriguing results have been recently proved, see
for instance Da Prato \& Debussche \cite{DD}, Mikulevicius \& Rozovski
\cite{MR}, Flandoli \& Romito \cite{FR} (see also \cite{FRcras}). We recall
here the framework constructed in \cite{FR} and prove some additional results.

We consider a viscous, incompressible, homogeneous, Newtonian fluid described
by the stochastic Navier-Stokes equations on the torus $\mathcal{T}=\left[
0,L\right]  ^{3}$, $L>0$,
\begin{equation}\label{e:stochNS}
\frac{\partial u}{\partial t}+\left(  u\cdot\nabla\right)  u+\nabla
p=\nu\triangle u+\sum_{i=1}^{\infty}\sigma_{i}h_{i}\left(  x\right)
\overset{\cdot}{\beta}_{i}\left(  t\right)
\end{equation}
with $\mathsf{div}\,u=0$ and periodic boundary conditions, with suitable
fields $h_{i}\left(  x\right)  $ and independent Brownian motions $\beta
_{i}\left(  t\right)  $. The 3D random vector field $u=u\left(  t,x\right)  $
is the velocity of the fluid and the random scalar field $p=p\left(
t,x\right)  $ is the pressure. To simplify the exposition, we avoid generality
and focus on one of the simplest set of assumptions:
$$
\sigma_{i}^{2}=\lambda_{i}^{-3}
$$
where $\lambda_{i}$ are the eigenvalues of the Stokes operator (see the next
section). This assumption also allows us to compare more closely the results
in Da Prato \& Debussche \cite{DD} and Flandoli \cite{F}. However, following
Flandoli \& Romito \cite{FR}, we could treat any power law for $\sigma_{i}$.
Under this assumption, one can associate a transition probability kernel
$P\left(  t,x,\cdot\right)  $ to equation \eqref{e:SNSeq}, which is the abstract
version of \eqref{e:stochNS}, in $D(A)$ (see the definitions in Section
\ref{ss:notations} below), satisfying the Chapman-Kolmogorov equation. In
other words, there exists a Markov selection in $D(A)$ for equation
\eqref{e:SNSeq}. To avoid misunderstandings, this does not mean that equation
\eqref{e:SNSeq} has been solved in $D(A)$ with continuous trajectories: this
would imply well posedness. What has been proved is that the law of weak
martingale solutions is supported on $D(A)$ for all times, with a number of
related additional properties, but a priori the typical trajectory may
sometimes blow-up in the topology of $D(A)$.

The transition probabilities $P\left(  t,x,\cdot\right)  $ are irreducible and
strong Feller, hence equivalent, in $D(A)$. These results and the existence of
$P\left(  t,x,\cdot\right)  $ have been proved first in Da Prato \& Debussche
\cite{DD} and Debussche \& Odasso \cite{DO} by a careful selection from the
Galerkin scheme. Then another proof by an abstract selection principle and the
local-in-time regularity of equation \eqref{e:SNSeq} has been given in Flandoli
\& Romito \cite{FR}. More precisely, first one proves the existence of a
Markov kernel $P\left(  t,x,\cdot\right)  $ by means of a general and abstract
method, then one proves that \textit{any} such kernel is irreducible and
strong Feller, hence equivalent, in $D(A)$.

We complement here the approach of \cite{FR} with two results. First, the
simple idea used in \cite{FR} to prove the strong Feller property is here
developed further, to show a weak form of Lipschitz continuity of $P\left(
t,x,\cdot\right)  $ in $x\in D(A)$. More precisely, we prove the estimate
\begin{equation}\label{e:logLipschitz1}
\left|  P\left(  t,x_{0}+h,\Gamma\right)  -P\left(  t,x_{0},\Gamma\right)
\right|  \leq\frac{C_{T}}{t\wedge1}(1+|Ax_{0}|^{6})|Ah|\log(|Ah|^{-1})
\end{equation}
for $t\in(0,T]$, $x_{0}$, $h\in D(A)$, with $\left|  Ah\right|  \leq1$. This
result has been proved in a stronger version in Da Prato \& Debussche
\cite{DD} for the transition kernel constructed from the Galerkin scheme, and
also in Flandoli \cite{F} for any Markov kernel associated to equation
\eqref{e:SNSeq}. In both cases the proof is based on the very powerful approach
introduced in \cite{DD} which however requires a considerable amount of
technical work. Here we give a rather elementary proof along the lines of
Flandoli \& Romito \cite{FR}, based on the following simple idea: given
$x_{0}$, $h\in D(A)$, for a short random time the solution is regular, unique
and differentiable in the initial conditions; then the propagation of
regularity in $x$ from small time to arbitrary time is due to the Markov
property. Unfortunately we cannot prove in this way the stronger estimate
obtained in \cite{DD} (where the right-hand-side of \eqref{e:logLipschitz1}
has the form $t^{-1+\varepsilon}(1+|Ax_{0}|^{2})|Ah|$), so our first result
here has mostly a pedagogical character, since the proof is conceptually very easy.

The second result, which follows from the same main estimates used to prove
\eqref{e:logLipschitz1}, is the equivalence
\[
P^{\left(  1\right)  }\left(  t,x,\cdot\right)  \sim P^{\left(  2\right)
}\left(  t^{\prime},x^{\prime},\cdot\right)
\]
for any $t$, $t^{\prime}>0$ and $x$, $x^{\prime}\in D(A)$, when $P^{\left(
i\right)  }\left(  t,x,\cdot\right)  $, $i=1,2$, are any two Markov transition
kernels associated to equation \eqref{e:SNSeq} in $D(A)$. We have not proved yet
the existence of invariant measures associated to such kernels\footnote{This
is apparently due to technical reasons and it is the subject of a work in
progress.}, but if we assume to have such invariant measures, it also follows
that they are equivalent. This result and the gradient estimates discussed
above could be steps to understand better the open question of well posedness
for equation \eqref{e:SNSeq}). In particular, it seems to be not so easy to
produce examples of stochastic differential equations without uniqueness but
where all Markov solutions are equivalent.

Among the open problems related to this research we mention the relation
between the regularity results for $P\left(  t,x,\cdot\right)  $ in the
initial condition discussed above and the properties of Malliavin derivatives,
investigated for stochastic 3D Navier-Stokes equations by Mikulevicius and
Rozovsky in \cite{MR0} and \cite{MR}.
%%
%%%%%%%%%%%%%%%%%%%%%%%%%%%%%%%%%%%%%%%%%%%%%%%%%%%%%%%%%%%%%%%%%%%%%%%%%%%%%%
%%
\section{Preliminaries}
\subsection{Notations}\label{ss:notations}
Denote by $\mathcal{T}=[0,1]^{3}$ the three-dimensional torus, and let
$\mathbb{L}^{2}\left(  \mathcal{T}\right)  $ be the space of vector fields
$u:\mathcal{T}\rightarrow\mathbb{R}^{3}$ with $L^{2}\left(  \mathcal{T}%
\right)  $-components. For every $\alpha>0$, let $\mathbb{H}^{\alpha}\left(
\mathcal{T}\right)  $ be the space of fields $u\in\mathbb{L}^{2}\left(
\mathcal{T}\right)  $ with components in the Sobolev space $H^{\alpha}\left(
\mathcal{T}\right)  =W^{\alpha,2}\left(  \mathcal{T}\right)  $.

Let $\mathcal{D}^{\infty}$ be the space of infinitely differentiable
divergence free periodic fields $u$ on $\mathcal{T}$, with zero mean. Let $H$
be the closure of $\mathcal{D}^{\infty}$ in the topology of $\mathbb{L}%
^{2}\left(  \mathcal{T}\right)  $: it is the space of all zero mean fields
$u\in\mathbb{L}^{2}\left(  \mathcal{T}\right)  $ such that $\mathsf{div}\,u=0$
and $u\cdot n$ on the boundary is periodic. We denote by $\left\langle
.,.\right\rangle _{H}$ and $\left|  .\right|  _{H}$ (or simply by
$\left\langle .,.\right\rangle $ and $\left|  .\right|  $) the usual $L^{2}%
$-inner product and norm in $H$. Let $V$ (resp. $D(A)$) be the closure of
$\mathcal{D}^{\infty}$ in the topology of $\mathbb{H}^{1}\left(
\mathcal{T}\right)  $ (in the topology of $\mathbb{H}^{2}\left(
\mathcal{T}\right)  $, respectively): it is the space of divergence free, zero
mean, periodic elements of $\mathbb{H}^{1}\left(  \mathcal{T}\right)  $
(respectively of $\mathbb{H}^{2}\left(  \mathcal{T}\right)  $). The spaces $V$
and $D(A)$ are dense and compactly embedded in $H$. From Poincar\'{e}
inequality we may endow $V$ with the norm $\left\|  u\right\|  _{V}^{2}%
:=\int_{\mathcal{T}}\left|  Du\left(  x\right)  \right|  ^{2}\,dx $.

Let $A:D(A)\subset H\rightarrow H$ be the operator $Au=-\triangle u$
(component wise). Since $A$ is a selfadjoint positive operator in $H$, there
is a complete orthonormal system $(h_{i})_{i\in\mathbb{N}}\subset H$ of
eigenfunctions of $A$, with eigenvalues $0<\lambda_{1}\leq\lambda_{2}%
\leq\ldots$ (that is, $Ah_{i}=\lambda_{i}h_{i}$). The fields $h_{i}$ in
equation \eqref{e:SNSeq} will be these eigenfunctions. We have
\[
\left\langle Au,u\right\rangle _{H}=\left\|  u\right\|  _{V}^{2}%
\]
for every $u\in D(A)$.

Let $V^{\prime}$ be the dual of $V$; with proper identifications we have
$V\subset H\subset V^{\prime}$ with continuous injections, and the scalar
product $\left\langle \cdot,\cdot\right\rangle _{H}$ extends to the dual
pairing $\left\langle \cdot,\cdot\right\rangle _{V,V^{\prime}}$ between $V $
and $V^{\prime}$. We may enlarge this scheme to $D(A)\subset V\subset H\subset
V^{\prime}\subset D(A)^{\prime}$. Let $B\left(  \cdot,\cdot\right)  :V\times
V\rightarrow V^{\prime}$ be the bi linear operator defined as
\[
\left\langle w,B\left(  u,v\right)  \right\rangle _{V,V^{\prime}}=\sum
_{i,j=1}^{3}\int_{\mathcal{T}}u_{i}\frac{\partial v_{j}}{\partial x_{i}}%
w_{j}\,dx
\]
for every $u,v,w\in V$. We shall repeatedly use the following inequality:
\begin{equation}\label{e:inequalityAB}
\left|  A^{1/2}B\left(  u,v\right)  \right|  _{H}\leq C_{0}\left|  Au\right|
\left|  Av\right|
\end{equation}
for $u,v\in D(A)$. The proof is elementary (see Flandoli \cite{FlaCime}).
\subsection{Definitions, assumptions and known results}
We (formally) rewrite equations \eqref{e:stochNS} as an abstract stochastic
evolution equation in $H$,
\begin{equation}\label{e:SNSeq}
du(t)+\left[  \nu Au(t)+B\left(  u(t),u(t)\right)  \right]  \, dt=\sum
_{i=1}^{\infty}\sigma_{i}h_{i}\,d\beta_{i}\left(  t\right)  .
\end{equation}
Let us set
$$
\Omega=C([0,\infty);D(A)')
$$
and denote by $\left(  \xi_{t}\right)  _{t\geq0}$ the canonical process on
$\Omega$, defined as $\xi_{t}\left(  \omega\right)  =\omega\left(  t\right)
$, by $\mathcal{F}$ the Borel $\sigma$-algebra in $\Omega$ and by
$\mathcal{F}_{t}$ the $\sigma$-algebra generated by the events $\left\{
\xi_{s}\in A\right\}  $ with $s\in\left[  0,t\right]  $ and $A$ a Borel set of
$D(A)^{\prime}$. Finally, denote by $\mathcal{B}(D(A))$ the Borel $\sigma
$-algebra of $D(A)$ and by $B_{b}(D(A))$ the set of all real valued bounded
measurable functions on $D(A))$.
\begin{definition}
Given a probability measure $\mu_{0}$ on $H$, we say that a probability
measure $P$ on $\left(  \Omega,\mathcal{F}\right)  $ is a solution to the
martingale problem associated to equation \eqref{e:SNSeq} with initial law
$\mu_{0}$ if
\begin{enumerate}
\item[\textbf{\footnotesize\textsf{\lbrack MP1\rbrack}}] $P[\xi\in L_{\text{loc}}^{\infty}([0,\infty);H)\cap L_{\text{loc}}^{2}([0,\infty);V)]=1$,
\item[\textbf{\footnotesize\textsf{\lbrack MP2\rbrack}}] for each $\varphi\in\mathcal{D}^{\infty}$ the process
$(M_{t}^{\varphi},\mathcal{F}_{t},P)_{t\geq0}$, defined $P$-a.\ s.\ on
$(\Omega,\mathcal{F})$ as
$$
M_{t}^{\varphi}:=\langle\xi_{t}-\xi_{0},\varphi\rangle_{H}+\int_{0}^{t}%
\nu\langle\xi_{s},A\varphi\rangle_{H}\,ds-\int_{0}^{t}\langle B(\xi
_{s},\varphi),\xi_{s}\rangle_{H}\,ds
$$
is a continuous square integrable martingale with quadratic variation
$$
\lbrack M^{\varphi}]_{t}=t\sum_{i\in\mathbb{N}}\sigma_{i}^{2}|\langle
\varphi,h_{i}\rangle|^{2},
$$
\item[\textbf{\footnotesize\textsf{\lbrack MP3\rbrack}}] the marginal of $P$ at time $0$ is $\mu_{0}$.
\end{enumerate}
\end{definition}
\begin{remark}\label{r:equivalence}
Among all test functions in property \textbf{\footnotesize\textsf{[MP2]}}, we can
choose $\varphi=h_{i}$. Set for all $i$, $\beta_{i}(t)=\frac{1}{\sigma_{i}%
}M_{t}^{h_{i}}$ (and $0$ if $\sigma_{i}=0$). The $(\beta_{i})_{i\in\mathbb{N}%
}$ are a sequence of independent standard Brownian motions. Under the
assumption $\sum_{i}\sigma_{i}^{2}<\infty$, the series $\sum_{i=1}^{\infty
}\sigma_{i}h_{i}\beta_{i}\left(  t\right)  $ defines an $H$-valued Brownian
motion on $\left(  \Omega,\mathcal{F},\mathcal{F}_{t},P\right)  $, that we
shall denote by $W\left(  t\right)  $. The canonical process $\left(  \xi
_{t}\right)  $ is a weak martingale solution of \eqref{e:SNSeq}, in the sense
that it satisfies \eqref{e:SNSeq} in the following weak form: there exists a
Borel set $\Omega_{0}\subset\Omega$ with $P\left(  \Omega_{0}\right)  =1$ such
that on $\Omega_{0}$ for every $\varphi\in\mathcal{D}^{\infty}$ and $t\geq0$
we have
\begin{equation}
\left\langle \xi_{t}-\xi_{0},\varphi\right\rangle _{H}+\int_{0}^{t}%
\nu\left\langle \xi_{s},A\varphi\right\rangle _{H}\,ds-\int_{0}^{t}%
\left\langle B\left(  \xi_{s},\varphi\right)  ,\xi_{s}\right\rangle
_{H}\,ds=\left\langle W\left(  t\right)  ,\varphi\right\rangle _{H}.
\end{equation}
\end{remark}
The following theorem is well known, see for instance the survey paper of
Flandoli \cite{FlaCime} and the reference therein.
\begin{theorem}\label{t:teoexist1}
Assume $\sum_{i}\sigma_{i}^{2}<\infty$. Let $\mu$ be a
probability measure on $H$ such that $\int_{H}\left|  x\right|  _{H}^{2}%
\mu\left(  dx\right)  <\infty$. Then there exists at least one solution to the
martingale problem with initial condition $\mu$.
\end{theorem}
\begin{definition}
We say that $P\left(  \cdot,\cdot,\cdot\right)  :[0,\infty)\times
D(A)\times\mathcal{B}\left(  D(A)\right)  \rightarrow\left[  0,1\right]  $ is
a Markov kernel in $D(A)$ of transition probabilities associated to equation
\eqref{e:stochNS} if $P\left(  \cdot,\cdot,\Gamma\right)  $ is Borel measurable
for every $\Gamma\in$ $\mathcal{B}\left(  D(A)\right)  $, $P\left(
t,x,\cdot\right)  $ is a probability measure on $\mathcal{B}\left(
D(A)\right)  $ for every $\left(  t,x\right)  \in\lbrack0,\infty)\times D(A)$,
the Chapman-Kolmogorov equation
\[
P\left(  t+s,x,\Gamma\right)  =\int_{D(A)}P\left(  t,x,dy\right)  P\left(
s,y,\Gamma\right)
\]
holds for every $t,s\geq0$, $x\in D(A)$, $\Gamma\in$ $\mathcal{B}\left(
D(A)\right)  $, and for every $x\in D(A)$ there is a solution $P_{x}$ on
$\left(  \Omega,F\right)  $ of the martingale problem associated to equation
\eqref{e:SNSeq} with initial condition $x$ such that
\[
P\left(  t,x,\Gamma\right)  =P_{x}\left[  \xi_{t}\in\Gamma\right]  \text{ for
all }t\geq0\text{.}
\]
\end{definition}
We recall the following result from Da Prato \& Debussche \cite{DD}, Debussche
\& Odasso \cite{DO} or Flandoli \& Romito \cite{FR}:
\begin{theorem}
There exists at least one Markov kernel $P\left(  t,x,\Gamma\right)  $ in
$D(A)$ of transition probabilities associated to equation \eqref{e:stochNS}.
\end{theorem}

We recall that a $P\left(  t,x,\Gamma\right)  $ is called \emph{irreducible}
in $D(A)$ if for every $t>0$, $x_{0}$, $x_{1}\in D(A)$, $\varepsilon>0$, we
have
\[
P\left(  t,x_{0},B_{A}\left(  x_{1},\varepsilon\right)  \right)  >0,
\]
where $B_{A}\left(  x_{1},\varepsilon\right)  $ is the ball in $D(A)$ of
centre $x_{1}$ and radius $\varepsilon$.

We say that $P\left(  t,x,\Gamma\right)  $ is \emph{strong Feller} in $D(A)$
if
\[
x\mapsto\int_{D(A)}\varphi\left(  y\right)  P\left(  t,x,dy\right)
\]
is continuous on $D(A)$ for every bounded measurable function $\varphi
:D(A)\rightarrow\mathbb{R}$ and for every $t>0$. It is well known (see for
example Da Prato \& Zabczyk \cite[Proposition 4.1.1]{DPZ96}) that
irreducibility and strong Feller in $D(A)$ imply that the laws $P\left(
t,x,\cdot\right)  $ are all mutually equivalent, as $\left(  t,x\right)  $
varies in $(0,\infty)\times D(A)$. Because of this equivalence property, we
say that $P\left(  t,x,\Gamma\right)  $ is \emph{regular}.

We recall also that $P\left(  t,x,\Gamma\right)  $ is called
\emph{stochastically continuous} in $D(A)$ if $\lim_{t\rightarrow0}P\left(
t,x,B_{A}\left(  x,\varepsilon\right)  \right)  =1$ for every $x\in D(A)$ and
$\varepsilon>0$.

In Da Prato \& Debussche \cite{DD}, the transition probability kernel
constructed by Galerkin approximations is proved to be stochastically
continuous, irreducible and strong Feller in $D(A)$, hence regular. More
generally (see Flandoli \& Romito \cite{FR}):
\begin{theorem}
Every Markov kernel $P\left(  t,x,\Gamma\right)  $ in $D(A)$ of transition
probabilities associated to equation \eqref{e:stochNS} is stochastically
continuous, irreducible and strong Feller in $D(A)$, hence regular.
\end{theorem}
%%
%%%%%%%%%%%%%%%%%%%%%%%%%%%%%%%%%%%%%%%%%%%%%%%%%%%%%%%%%%%%%%%%%%%%%%%%%%%%
%%
\section{The Log-Lipschitz estimate}
\begin{theorem}\label{t:loglipschitz}
Let $P\left(  t,x,\Gamma\right)  $ be a Markov kernel
in $D(A)$ of transition probabilities associated to equation \eqref{e:stochNS}.
Then, given $T>0$, there is a constant $C_{T}$ such that the inequality
\[
\left|  P\left(  t,x_{0}+h,\Gamma\right)  -P\left(  t,x_{0},\Gamma\right)
\right|  \leq\frac{C_{T}}{t\wedge1}(1+|Ax_{0}|^{6})|Ah|\log(|Ah|^{-1})
\]
holds for every $t\in(0,T]$, $x_{0}$, $h\in D(A)$, with $|Ah|\leq1$, and
$\Gamma\in\mathcal{B}\left(  D(A)\right)  $.
\end{theorem}

We explain here only the logical skeleton of the proof, which is very simple.
The two main technical ingredients will be treated in the next two separate
subsections. The first idea is to decompose:
\begin{align*}
&  P\left(  t,x_{0}+h,\Gamma\right)  -P\left(  t,x_{0},\Gamma\right)  =\\
&  \qquad=\int_{D(A)}\left[  P\left(  \varepsilon,x_{0}+h,dy\right)  -P\left(
\varepsilon,x_{0},dy\right)  \right]  P\left(  t-\varepsilon,y,\Gamma\right)
.
\end{align*}
To shorten some notation, let us write
\[
\left(  P_{t}\varphi\right)  \left(  x\right)  =\int_{D(A)}\varphi\left(
y\right)  P\left(  t,x,dy\right)
\]
so, with the function $\varphi\left(  x\right)  =\mathbf{1}_{\{x\in\Gamma\}}$
the previous identity reads
\begin{align}\label{e:decomposition}
&  \left(  P_{t}\varphi\right)  \left(  x_{0}+h\right)  -\left(  P_{t}%
\varphi\right)  \left(  x_{0}\right)  =\left(  P_{\varepsilon}\left(
P_{t-\varepsilon}\varphi\right)  \right)  \left(  x_{0}+h\right)  -\left(
P_{\varepsilon}\left(  P_{t-\varepsilon}\varphi\right)  \right)  \left(
x_{0}\right)  .
\end{align}
It is now sufficient to estimate
\[
\left(  P_{\varepsilon}\psi\right)  \left(  x_{0}+h\right)  -\left(
P_{\varepsilon}\psi\right)  \left(  x_{0}\right)
\]
uniformly in $\psi\in B_{b}\left(  D(A)\right)  $. The value of $\varepsilon$
has to be chosen depending on the size of $x_{0}$ and $h$, as we shall see.

The second idea is to use an \emph{initial coupling}: we introduce the
equation with cut-off $\chi_{R}( \left|  Au\right|  ^{2})$, where $\chi
_{R}\left(  r\right)  :\left[  0,\infty\right)  \rightarrow\left[  0,1\right]
$ is a non-increasing smooth function equal to 1 over $\left[  0,R\right]  $,
to 0 over $\left[  R+2,\infty\right)  $, and with derivative bounded by 1. The
equation is
\begin{equation}\label{e:stochNSregular}%
\begin{array}
[c]{ll}%
du+\left[  Au+B(u,u)\chi_{R}\left(  \left|  Au\right|  ^{2}\right)  \right]
\,dt =\sum_{i=1}^{\infty}\sigma_{i}h_{i}\,d\beta_{i}\left(  t\right)  , & \\
u\left(  0\right)  =x. &
\end{array}
\end{equation}
The definition of martingale problem for this equation is the same (with
obvious adaptations) as the definition given above for equation \eqref{e:stochNS}.
Let $\tau_{R}:\Omega\rightarrow\left[  0,\infty\right]  $ be defined as
\[
\tau_{R}\left(  \omega\right)  =\inf\left\{  t\geq0:\left|  A\omega\left(
t\right)  \right|  \geq R\right\}  .
\]
We recall the following result from Flandoli \& Romito \cite[Lemma 5.11]{FR}:

\begin{lemma}
For every $x\in D(A)$ there is a unique solution $P_{x}^{\left(  R\right)  }$
of the martingale problem associated\ to equation \eqref{e:stochNSregular},
with the additional property
\[
P_{x}^{\left(  R\right)  }\left[  \xi\in C\left(  \left[  0,\infty\right)
;D(A)\right)  \right]  =1.
\]
Let $P_{x}$ be any solution on $\left(  \Omega,F\right)  $ of the martingale
problem associated to equation \eqref{e:SNSeq} with initial condition $x$. Then
\[
\mathbb{E}^{P_{x}^{\left(  R\right)  }}\left[  \varphi\left(  \xi_{t}\right)
\mathbf{1}_{\left\{  \tau_{R}\geq t\right\}  }\right]  =\mathbb{E}^{P_{x}%
}\left[  \varphi\left(  \xi_{t}\right)  \mathbf{1}_{\left\{  \tau_{R}\geq
t\right\}  }\right]
\]
for every $t\geq0$ and $\varphi\in B_{b}\left(  D(A)\right)  $.
\end{lemma}

Introduce the notation
\[
(P_{t}^{(R)}\varphi) \left(  x\right)  =\mathbb{E}^{P_{x}^{\left(  R\right)
}}\left[  \varphi\left(  \xi_{t}\right)  \right]  .
\]
The previous lemma implies that for every $\psi\in B_{b}\left(  D(A)\right)  $
we have
\begin{equation}\label{e:confronto}
|(P_{\varepsilon}\psi)(x)-(P_{\varepsilon}^{(R)}\psi)(x)| \leq2P_{x}[\tau
_{R}<\varepsilon]\,\|\psi\|_{\infty}. 
\end{equation}
Summarising:
\begin{corollary}
For every $x_{0}$, $h\in D(A)$ and $\psi\in B_{b}\left(  D(A)\right)  $ we
have
\begin{align*}
|(P_{\varepsilon}\psi)(x_{0}+h)-(P_{\varepsilon}\psi)(x_{0})| &  \leq2\left(
P_{x_{0}+h}[\tau_{R}<\varepsilon]+P_{x_{0}}[\tau_{R}<\varepsilon]\right)
\Vert\psi\Vert_{\infty}\\
&  \quad+\left|  (P_{\varepsilon}^{(R)}\psi)(x_{0}+h)-(P_{\varepsilon}^{(R)}\psi)(x_{0})\right|.
\end{align*}
\end{corollary}
Let us give now the proof of Theorem \ref{t:loglipschitz}. Assume $t\in(0,T]$,
$x_{0}$, $h\in D(A)$ be given, with $|Ah|\leq1$. Let $K>0$ be such that $|A
x_{0}|+1\leq K$. We have $|A(x_{0}+h)|\leq K$, so we may apply Proposition
\ref{p:propblowuptime} below to both $x_{0}$ and $x_{0}+h$. We thus get, for
$\varepsilon\in(0,\frac{1}{5C^{*}K^{2}})$, where $C^{\ast}>0$ is the constant
defined by \eqref{e:defCstar}, we have
$$
P_{x_{0}+h}\left[  \tau_{2K}<\varepsilon\right]  +P_{x_{0}}\left[  \tau
_{2K}<\varepsilon\right]  \leq2C_{\#}\mathrm{e}^{-\eta_{\#}\frac{K^{2}%
}{4\varepsilon}}.
$$
Given $h$, $K$ and $t$ as above, let us look for a value $\varepsilon
\in(0,\frac{1}{5C^{*}K^{2}})$ such that $\varepsilon\le t$ and the latter
exponential quantity is smaller than $|Ah|$. We impose
$$
\eta_{\#}\frac{K^{2}}{4\varepsilon}\geq\log(|Ah|^{-1})
$$
hence it is sufficient to take
\begin{equation}\label{e:epsilonbound}
\varepsilon\le\frac{\eta_{\#}K^2}{4\log(|Ah|^{-1})}\wedge\frac{t}{2}\wedge\frac{1}{5C^{*}K^{2}}.
\end{equation}
We have proved so far the first claim of the following lemma. The second claim
is a simple consequence of \eqref{e:decomposition} and the previous corollary.
\begin{lemma}
Given $t>0$, $x_{0}$, $h\in D(A)$, with $|Ah|\le1$, and $\Gamma\in\mathcal{B}(D(A))$,
if $\varepsilon$ is chosen as in \eqref{e:epsilonbound}, then
$$
P_{x_{0}+h}\left[  \tau_{2K}<\varepsilon\right]  +P_{x_{0}}\left[  \tau
_{2K}<\varepsilon\right]  \leq{2C_{\#}}|Ah|
$$
and for $\varphi(x)=\mathbf{1}_{\{x\in\Gamma\}}$ and $\psi=P_{t-\varepsilon}\varphi$,
$$
|P_{t}\varphi(x_{0}+h)-P_{t}\varphi(x_{0})|\leq4C_{\#}|Ah|\,\Vert\varphi
\Vert_{\infty}+|P_{\varepsilon}^{(2K)}\psi(x_{0}+h)-P_{\varepsilon}^{(2K)}%
\psi(x_{0})|.
$$
\end{lemma}

Finally, from Proposition \ref{p:propderivative} below, renaming the constant
$C$, with $\varphi\left(  x\right)  =\mathbf{1}_{\{x\in\Gamma\}}$ and
$\psi=P_{t-\varepsilon}\varphi$,
\[
\left|  ( P_{\varepsilon}^{\left(  2K\right)  }\psi) \left(  x_{0}+h\right)
-( P_{\varepsilon}^{(2K) }\psi) \left(  x_{0}\right)  \right|  \leq\frac
{C}{\varepsilon}\left|  Ah\right|  e^{CK^{6}\varepsilon}.
\]
Thus, for $\varepsilon$ as in \eqref{e:epsilonbound}, we get
\[
\left|  P_{t}\varphi\left(  x_{0}+h\right)  -P_{t}\varphi\left(  x_{0}\right)
\right|  \leq4C_{\#}|Ah|+\frac{C}{\varepsilon}|Ah|\mathrm{e}^{CK^{6}%
\varepsilon}.
\]
Let us further restrict ourselves to
\[
\varepsilon\leq\frac{\eta_{\#}K^{2}}{4\log(|Ah|^{-1})}\wedge\frac{t}{2}
\wedge\frac{1}{5C^{*}K^{2}}\wedge\frac{1}{K^{6}},
\]
so that we have
\[
|P_{t}\varphi(x_{0}+h)-P_{t}\varphi(x_{0})| \leq4C_{\#}|Ah|+\frac
{C}{\varepsilon}|Ah|.
\]
The choice
\[
\varepsilon=C\frac{t\wedge1}{K^{6}\log(|Ah|^{-1})}
\]
is admissible for a suitable constant $C>0$, and we finally get \eqref{e:logLipschitz1}.
The proof of Theorem \ref{t:loglipschitz} is complete.
\subsection{Probability of blow-up}
\begin{proposition}\label{p:propblowuptime}
Let $K\geq1$ and assume that $x_{0}\in D(A)$ and
$\varepsilon>0$ are given such that $|Ax_{0}|\leq K$ and $\varepsilon\leq
\frac{1}{5C^{\ast}K^{2}}$, where $C^{\ast}$ is the constant defined in
\eqref{e:defCstar}. Then
$$
P_{x_{0}}[\tau_{2K}<\varepsilon]\le C_{\#}\mathrm{e}^{-\eta_{\#}\frac{K^{2}}{4\varepsilon}},
$$
for suitable universal constants $\eta_{\#}>0$ and $C_{\#}>0$.
\end{proposition}
\begin{proof}
From Corollary \ref{c:corollaryappendix} we know that if $\varepsilon\leq
\frac{1}{5C^{\ast}K^{2}}$ and $|Ax_{0}|\leq K$, then one has
\[
\theta_{\varepsilon}^{2}\leq\frac{1}{4}K^{2}\quad\Rightarrow\quad\left|
Au(s)\right|  <2K\text{ for }s\in\left[  0,\varepsilon\right]  \quad
\Rightarrow\quad\tau_{2K}\geq\varepsilon,
\]
where $\theta_{\varepsilon}$ is defined in Section \eqref{ss:deftheta}.
Therefore, with the constraints $|Ax_0|\le K$ and
$\varepsilon\le\frac{1}{5C^{\ast}K^{2}}$, by Proposition \ref{p:Ztails}
one gets
$$
P_{x_{0}}\left[  \tau_{2K}<\varepsilon\right]  \leq P_{x_{0}}\left[
\Theta_{\varepsilon}^{2}>\frac{1}{4}K^{2}\right]  \leq C_{\#}\mathrm{e}%
^{-\eta_{\#}\frac{K^{2}}{4\varepsilon}}.
$$
\end{proof}
\subsection{Derivative of the regularised problem}
Here we show the regularity of the transition semigroup associated to the
regularised problem \eqref{e:stochNSregular}.
\begin{proposition}\label{p:propderivative}
For every $R\geq1$ and $x_{0}$, $h\in D(A)$,
$$
|(P_\varepsilon^{(R)}\psi)(x_0+h)-(P_\varepsilon^{(R)}\psi)(x_0)|
\le\frac{C\|\psi\|_\infty}{\varepsilon}|Ah|\mathrm{e}^{CR^6\varepsilon},
$$
where $C$ is a universal constant.
\end{proposition}
\begin{proof}
We write the following computations for the limit problem but the
understanding is that we do it on the Galerkin approximations. For every
$\psi\in B_b(H)$, $\varepsilon>0$, from the
Bismut-Elworthy-Li formula (see Da Prato \& Zabczyk \cite{DPZ96}),
\begin{multline*}
|(P_\varepsilon^{(R)}\psi)(x_0+h)-(P_\varepsilon^{(R)}\psi)(x_0)|\le\\
\le\frac{C\|\psi\|_\infty}{\varepsilon}\sup_{\eta\in[0,1]}\mathbb{E}\bigl[\bigl(\int_0^\varepsilon|A^{\frac32}D_h u_{x_0+\eta h}^{(R)}(s)|^2\,ds\bigr)^{\frac12}\bigr],
\end{multline*}
where, for each $R\geq1$ and $x\in D(A)$, $u_{x}^{(R)}$ is the solution,
starting at $x$, of problem \eqref{e:stochNSregular}. From the regularised
equation we have
\begin{align*}
&\frac12\frac{d}{dt}|AD_h u_x^{(R)}(t)|^2+|A^{\frac32}D_h u_x^{(R)}(t)|^2\le\\
&\qquad\le    \chi_R(|Au_x^{(R)}(t)|^2)|\langle AD_h u_x^{(R)},AB(D_h u_x^{(R)},u_x^{(R)})+AB(u_x^{(R)},D_h u_x^{(R)})\rangle|\\
&\qquad\quad +2\chi_R'(|Au_x^{(R)}(t)|^2)\langle Au_x^{(R)},AD_h u_x^{(R)}\rangle|\langle AD_h u_x^{(R)},AB(u_x^{(R)},u_x^{(R)})\rangle|\\
&\qquad\le    C\chi_R(|Au_x^{(R)}(t)|^2)|A^{\frac32}D_h u_x^{(R)}(t)|\,|AD_h u_x^{(R)}(t)|\,|Au_x^{(R)}(t)|\\
&\qquad\quad +C\chi_R'(|Au_x^{(R)}(t)|^2)|Au_x^{(R)}(t)|^3|AD_h u_x^{(R)}(t)|\,|A^{\frac32}D_h u_x^{(R)}(t)|\\
&\qquad\le    \frac12|A^{\frac32}D_h u_x^{(R)}(t)|^2+C\chi_R^2(|Au_x^{(R)}(t)|^2)|AD_h u_x^{(R)}(t)|^2|Au_x^{(R)}(t)|^2\\
&\qquad\quad +C\chi_R'(|Au_x^{(R)}(t)|^2)^2|AD_h u_x^{(R)}(t)|^2|Au_x^{(R)}(t)|^6\\
&\qquad\le    \frac12|A^{\frac32}D_h u_x^{(R)}(t)|^2+CR^6|AD_h u_x^{(R)}(t)|^2.
\end{align*}
Thus
$$
\frac{1}{2}\frac{d}{dt}\left|  AD_{h}u_{x}^{(R)}(t)\right|  ^{2}+\frac{1}%
{2}\left|  A^{\frac{3}{2}}D_{h}u_{x}^{(R)}(t)\right|  ^{2}\leq CR^{6}\left|
AD_{h}u_{x}^{(R)}(t)\right|  ^{2}.
$$
This implies
$$
\left|  AD_{h}u_{x}^{(R)}(t)\right|  ^{2}\leq\mathrm{e}^{CR^{6}t}\left|
Ah\right|  ^{2}%
$$
and
$$
\int_{0}^{\varepsilon}\left|  A^{\frac{3}{2}}D_{h}u_{x_{0}+\eta h}%
^{(R)}(s)\right|  ^{2}\,ds\leq\left|  Ah\right|  ^{2}\left(  1+\int
_{0}^{\varepsilon}CR^{6}\mathrm{e}^{CR^{6}s}\,ds\right)  =\left|  Ah\right|
^{2}\mathrm{e}^{CR^{6}\varepsilon}.
$$
Thus
$$
\left|  \left(  P_{\varepsilon}^{\left(  R\right)  }\psi\right)  \left(
x_{0}+h\right)  -\left(  P_{\varepsilon}^{\left(  R\right)  }\psi\right)
\left(  x_{0}\right)  \right|  \leq\frac{C\left\|  \psi\right\|  _{\infty}%
}{\varepsilon}\left|  Ah\right|  \mathrm{e}^{CR^{6}\varepsilon}.
$$
The proposition is proved.
\end{proof}

\section{Equivalence of all transition probabilities}

To make the following statement independent of previous results, we shall
assume stochastic continuity, irreducibility and the strong Feller property in
the theorem below, but we recall that these properties have been proved for
every Markov kernel in $D(A)$ associated to equation \eqref{e:stochNS}, under
the assumptions of the introduction.

\begin{theorem}
Let $P^{\left(  i\right)  }\left(  t,x,\Gamma\right)  $ be two Markov kernels
in $D(A)$ of transition probabilities associated to equation \eqref{e:stochNS}.
Assume they are stochastically continuous, irreducible and strong Feller in
$D(A)$. Then the probability measures $P^{\left(  1\right)  }\left(
t,x,\cdot\right)  $ and $P^{\left(  2\right)  }\left(  t^{\prime},x^{\prime
},\cdot\right)  $ are equivalent, for any $t$, $t^{\prime}>0$ and $x$,
$x^{\prime}\in D(A)$.
\end{theorem}

\begin{proof}
\underline{\emph{Step 1}.} Let $\Gamma$ be a Borel set in $D(A)$ such that
$P^{\left(  2\right)  }\left(  t_{0},x_{0},\Gamma\right)  =0$ for some
$t_{0}>0$, $x_{0}\in D(A)$. It is sufficient to prove that $P^{\left(
1\right)  }\left(  t_{0},x_{0},\Gamma\right)  =0$. We know that $P^{\left(
2\right)  }\left(  t,x,\Gamma\right)  =0$ for every $t>0$, $x\in D(A)$.

\noindent\underline{\emph{Step 2}.} Since both $P^{(1)}(\cdot,\cdot,\cdot)$
and $P^{(2)}(\cdot,\cdot,\cdot)$ satisfy \eqref{e:confronto},
\[
P^{(1)}(t,x,\Gamma)=|P^{(1)}(t,x,\Gamma)-P^{(2)}(t,x,\Gamma)|\leq2(P_{x}%
^{(1)}[\tau_{R}<t]+P_{x}^{(2)}[\tau_{R}<t]).
\]
Now, for every pair $(\varepsilon,x)$, with $\varepsilon>0$ and $x\in D(A)$,
such that $5C^{\ast}(1+|Ax|)^{2}\varepsilon\leq1$ (the constant $C^{\ast}$ is
defined in \eqref{e:defCstar}, in the appendix), Proposition \ref{p:propblowuptime}
implies that
$$
P^{(1)}(\varepsilon,x,\Gamma)\leq2C_{\#}\mathrm{e}^{-\eta_{\#}\frac
{(1+|Ax|)^{2}}{4\varepsilon}}\leq2C_{\#}\mathrm{e}^{-\frac{1}{4\varepsilon
}\eta_{\#}}.
$$

\noindent\underline{\emph{Step 3}.} For every $\varepsilon<\frac{1}{5C^{\ast}%
}$, set $A_{\varepsilon}=\{x\in D(A):5C^{\ast}(1+|Ax|)^{2}\varepsilon\leq1\}$,
then by the Markov property and the previous step,
\begin{align*}
P^{(1)}(t_{0}+\varepsilon,x_{0},\Gamma) &  =\int_{A_{\varepsilon}^{c}}%
P^{(1)}(\varepsilon,x,\Gamma)\,P^{(1)}(t_{0},x_{0},dx)\\
&  \quad+\int_{A_{\varepsilon}}P^{(1)}(\varepsilon,x,\Gamma)\,P^{(1)}%
(t_{0},x_{0},dx)\\
&  \leq2C_{\#}\mathrm{e}^{-\frac{1}{4\varepsilon}\eta_{\#}}+P^{(1)}%
(t_{0},x_{0},A_{\varepsilon}^{c})
\end{align*}
Since $P^{\left(  1\right)  }\left(  s,x_{0},D(A)\right)  =1$, we have
$P^{(1)}(t_{0},x_{0},A_{\varepsilon}^{c})\longrightarrow0$, as $\varepsilon
\rightarrow0$, and thus
\[
\lim_{\varepsilon\rightarrow0}P^{\left(  1\right)  }\left(  t_{0}%
+\varepsilon,x_{0},\Gamma\right)  =0.
\]

\noindent\underline{\emph{Step 4}.} By the Markov property, for every
neighborhood $G$ of $x_{0}$ in $D(A)$,
\begin{align*}
P^{(1)}(t_{0}+\varepsilon,x_{0},\Gamma) &  =\int P^{(1)}(t_{0},y,\Gamma
)\,P^{(1)}(\varepsilon,x_{0},dy)\\
&  \geq P^{(1)}(\varepsilon,x_{0},G)\inf_{y\in G}P^{(1)}(t_{0},y,\Gamma).
\end{align*}
Since the kernel $P^{(1)}$ is stochastically continuous, $P^{(1)}(\varepsilon,x_{0},G)$
converges to $1$, as $\varepsilon\rightarrow0$, and so, by the previous step,
$\inf_{y\in G}P^{(1)}(t_{0},y,\Gamma)\longrightarrow0$ as $\varepsilon\rightarrow0$.
By the strong Feller property, the map $y\mapsto P^{(1)}(t_{0},y,\Gamma)$ is
continuous, hence in conclusion $P^{(1)}(t_{0},x_{0},\Gamma)=0$. The proof
is complete.
\end{proof}
%%
%%%%%%%%%%%%%%%%%%%%%%%%%%%%%%%%%%%%%%%%%%%%%%%%%%%%%%%%%%%%%%%%%%%%%%%%%%%%%
%%
\section{Conclusion and remarks}
We have proved that the transition probabilities associated to any Markov
selection are all equivalent to each other. However, the problem of uniqueness
of Markov selections remains open. We stress that it would imply uniqueness of
solutions to the martingale problem, by the argument that one can find in
Stroock \& Varadhan \cite[Theorem 12.2.4]{SV}.

The estimates proved in this work allows us at least to state a sufficient
condition for uniqueness of Markov selections. The proof is inspired to a well
known proof in semigroup theory as well as to the proof of uniqueness given by
Bressan and co-authors (see for instance \cite{Bre}).
\begin{proposition}
Assume that a Markov selection $\left(  P_{x}\right)  _{x\in D(A)}$ has the
following property: for every $t>0$ and $x\in D(A)$,
\[
\lim_{n\rightarrow\infty}\sum_{k=1}^{n}P\left(  t-\frac{k}{n}t,x,B_{A}\left(
0,\sqrt{\frac{n}{t}}\right)  ^{c}\right)  =0
\]
where $B_{A}\left(  0,n\right)  $ is the ball in $D(A)$ of radius $n$. Then
$\left(  P_{x}\right)  _{x\in D(A)}$ coincides with any other Markov selection.
\end{proposition}
\begin{proof}
Let $\left(  Q_{x}\right)  _{x\in D(A)}$ be another Markov selection. Let us
rewrite, for $\varphi\in C_{b}\left(  D(A)\right)  $:
\begin{align*}
P_{t}\varphi-Q_{t}\varphi &  =P_{t-\frac{t}{n}}P_{\frac{t}{n}}\varphi
-P_{t-\frac{t}{n}}Q_{\frac{t}{n}}\varphi\\
&  +P_{t-\frac{t}{n}}Q_{\frac{t}{n}}\varphi-Q_{t-\frac{t}{n}}Q_{\frac{t}{n}%
}\varphi
\end{align*}
and so on iteratively until we have
$$
P_{t}\varphi-Q_{t}\varphi=\sum_{k=1}^{n}P_{t-\frac{kt}{n}}\left(  P_{\frac
{t}{n}}\psi_{\frac{\left(  k-1\right)  t}{n}}-Q_{\frac{t}{n}}\psi
_{\frac{\left(  k-1\right)  t}{n}}\right)
$$
where $\psi_{s}=Q_{s}\varphi$. We have, by using \eqref{e:confronto} and
Proposition \ref{p:propblowuptime},
\begin{align*}
&|P_{t-\frac{kt}{n}}(P_{\frac{t}{n}}\psi_{\frac{(k-1)t}{n}}-Q_{\frac{t}{n}}\psi_{\frac{(k-1)t}{n}})(x)|=\\
&\qquad  =   \bigl|\mathbb{E}^{P_x}\bigl[(P_{\frac{t}{n}}\psi_{\frac{(k-1)t}{n}}-Q_{\frac{t}{n}}\psi_{\frac{(k-1)t}{n}})(\xi_{t-\frac{kt}{n}})\bigl]\bigr|\\
&\qquad\le   \mathbb{E}^{P_x}\bigl[\bigl|(P_{\frac{t}{n}}\psi_{\frac{(k-1)t}{n}}-Q_{\frac{t}{n}}\psi_{\frac{(k-1)t}{n}})(\xi_{t-\frac{kt}{n}})\bigr|\mathbf{1}_{\{\xi_{t-\frac{k}{n}t}\in A_{\frac{t}{n}}\}}\bigr]\\
&\qquad\quad+\mathbb{E}^{P_x}\bigl[\bigl|(P_{\frac{t}{n}}\psi_{\frac{(k-1)t}{n}}-Q_{\frac{t}{n}}\psi_{\frac{(k-1)t}{n}})(\xi_{t-\frac{kt}{n}})\bigr|\mathbf{1}_{\{\xi_{t-\frac{k}{n}t}\in A_{\frac{t}{n}}^{c}\}}\bigr]\\
&\qquad\le   4C_{\#}\mathrm{e}^{-\frac{n}{t}\eta_{\#}}+2P_x[\xi_{t-\frac{k}{n}t}\in A_{\frac{t}{n}}^c]\\
&\qquad\le   4C_{\#}\mathrm{e}^{-\frac{n}{t}\eta_{\#}}+2P(t-\frac{k}{n}t,x,B_A(0,\sqrt{\frac{n}{t}})^c),
\end{align*}
where $A_t=\{5C^*t(1+|Ax|)^2\le1\}$ and, roughly,
$A_{\frac{t}{n}}\approx B_A(0,\sqrt{\frac{n}{t}})$. Hence
$$
|P_t\varphi(x)-Q_t\varphi(x)|
\le 4n\,C_{\#}\mathrm{e}^{-\frac{n}{t}\eta_{\#}}
   +2\sum_{k=1}^{n}P(t-\frac{k}{n}t,x,B_A(0,\sqrt{\frac{n}{t}})^c)
$$
which completes the proof of the proposition.
\end{proof}

The criterion of this proposition is apparently not really useful at the
present stage of our understanding. Indeed, if we apply Chebichev inequality
we get the sufficient condition
\[
\lim_{n\rightarrow\infty}\sum_{k=1}^{n}\left(  \frac{t}{n}\right)
^{1+\varepsilon}\mathbb{E}^{P_{x}}\left[  \left|  A\xi_{t-\frac{kt}{n}%
}\right|  ^{2\left(  1+\varepsilon\right)  }\right]  =0
\]
with is implied by the condition
\[
\mathbb{E}^{P_{x}}\left[  \int_{0}^{t}\left|  A\xi_{s}\right|  ^{2\left(
1+\varepsilon\right)  }\,ds\right]  <\infty
\]
which however would easily imply the well posedness of the 3D Navier-Stokes
equation by direct estimates of the difference of two solutions.
%%
%%%%%%%%%%%%%%%%%%%%%%%%%%%%%%%%%%%%%%%%%%%%%%%%%%%%%%%%%%%%%%%%%%%%%%%%%%%%%%%
%%
\appendix
\section{Appendix}
\subsection{A exponential tail estimate for the Stokes problem}\label{ss:deftheta}
Consider the following Stokes problem
$$
dZ+AZ\,dt=A^{-\frac32}\,dW,
\qquad
Z(0)=0,
$$
and set $\Theta_{t}=\sup_{s\in[0,t]}|AZ(s)|$. The next result is well known,
but we give a proof to keep track of the dependence on the constants of
interest in this paper.
\begin{proposition}\label{p:Ztails}
There exist $\eta_{\#}>0$ and $C_{\#}>0$ such that for every
$K\geq\frac{1}{2}$ and $\varepsilon>0$,
$$
\mathbb{P}\bigl[\Theta_{\varepsilon}\geq K]\leq C_{\#}\mathrm{e}^{-\eta
_{\#}\frac{K^{2}}{\varepsilon}}.
$$
\end{proposition}
\begin{proof}
\underline{\emph{Step 1}.} Set $y(t)=\varepsilon^{-\frac{1}{2}}Z(\varepsilon
t)$, then it is easy to see that $y$ solves the equation $dy+\varepsilon
Ay\,dt=Q^{\frac{1}{2}}\,dW$. Next, fix a value $\alpha\in(\frac{1}{6},\frac
{1}{4})$, then by the \emph{factorisation method} (see Da Prato \& Zabczyk
\cite[Chapter 5]{DPZ92}),
\[
y(t)=\int_{0}^{t}\mathrm{e}^{-\varepsilon(t-s)A}\,dW_{s}=C_{\alpha}\int
_{0}^{t}\mathrm{e}^{-\varepsilon(t-s)A}(t-s)^{\alpha-1}Y(s)\,ds,
\]
where $Y(s)=\int_{0}^{s}\mathrm{e}^{-\varepsilon(s-r)A}(s-r)^{-\alpha}%
\,dW_{r}$ and $C_{\alpha}$ denotes a generic constant depending only on
$\alpha$ (it will keep changing value along the proof). For every $t\in(0,1]$,
since $\alpha>\frac{1}{6}$, it follows from H\"{o}lder's inequality that
\[
|Ay(t)|_{H}\leq C_{\alpha}\int_{0}^{t}(t-s)^{\alpha-1}|AY(s)|_{H}\,ds\leq
C_{\alpha}\left(  \int_{0}^{1}|AY(s)|_{H}^{6}\,ds\right)  ^{\frac{1}{6}}.
\]
In conclusion, since $\varepsilon^{-1}\Theta_{\varepsilon}^{2}=\sup
_{t\in\lbrack0,1]}|Ay(t)|_{H}^{2}$, it follows by the above inequality and
standard arguments that
\begin{equation}\label{e:Zchebi}
\mathbb{P}[\Theta_{\varepsilon}\ge K]
\le\mathrm{e}^{-\frac{a}{\varepsilon}K^2}\mathbb{E}\Bigl[\exp\bigl(\tilde{a}(\int_0^1|AY(s)|^6\,ds)^{\frac13}\bigr)\Bigr],
\end{equation}
with a constant $a$ that will be specified later (and $\tilde{a}=a\,C_{\alpha}$).

\noindent\underline{\emph{Step 2}.} In order to estimate the expectation in
\eqref{e:Zchebi}, notice that
\begin{align}\label{e:Zserie}
&\exp\bigl(\tilde{a}\bigl[\int_0^1|AY(s)|^6\,ds\bigr]^{\frac13}\bigr)
=\sum_{n=0}^\infty\frac{{\tilde{a}}^n}{n!}\bigl[\int_0^1|AY(s)|_H^6\,ds\bigr]^{\frac{n}3}\leq\notag\\
&\qquad\le\sum_{n=0}^2\frac{{\tilde{a}}^n}{n!}\bigl[\int_0^1|AY(s)|_H^6\bigr]^{\frac{n}3}
         +\int_0^1\sum_{n=3}^\infty\frac{{\tilde{a}}^n}{n!}|AY(s)|_H^{2n}\\
&\qquad\le\tilde{a}\bigl[\int_0^1|AY(s)|^{10}\,ds\bigr]^{\frac15}
         +\frac{\tilde{a}^2}2\bigl[\int_0^1|AY(s)|^8\,ds\bigr]^{\frac12}
         +\int_0^1\mathrm{e}^{\tilde{a}|AY(s)|_H^2}\,ds\notag
\end{align}

\noindent\underline{\emph{Step 3}.} Now, $AY(s)$ is a centered Gaussian
process with covariance (cfr. proof of Theorem 5.9 in Da Prato \& Zabczyk
\cite{DPZ92})
\[
\tilde{Q}_{s}=\int_{0}^{s}(s-r)^{-2\alpha}A^{-1}\mathrm{e}^{-2\varepsilon
(s-r)A}\,dr,
\]
so that, by Proposition 2.16 of \cite{DPZ92},
\[
\mathbb{E}[\mathrm{e}^{\tilde{a}|AY(s)|_{H}^{2}}]=\mathrm{e}^{-\frac{1}%
{2}\mathrm{Tr}[\log(1-2\tilde{a}\tilde{Q}_{s})]},
\]
provided that $\tilde{a}\leq\inf_{\lambda\in\sigma(\tilde{Q}_{s})}\frac
{1}{2\lambda}$, where $\sigma(\tilde{Q}_{s})$ is the spectrum of $\tilde
{Q}_{s}$. Similarly, $\mathbb{E}|AY(s)|^{2p}=C_{p}(\mathrm{Tr}(\tilde{Q}%
_{s}))^{p}$, for all integers $p$.

In order to choose a suitable value of $a$, let $\mu\in\sigma(\tilde{Q}_{s})$,
then there is a eigenvalue $\lambda$ of $A$ such that $\mu=\mu(\lambda)$ is
given by
\[
\mu=\lambda^{-1}\int_{0}^{s}r^{-2\alpha}\mathrm{e}^{-2r\lambda\varepsilon
}\,dr=\lambda^{-2+2\alpha}(2\varepsilon)^{-(1-2\alpha)}\int_{0}^{2\lambda
\varepsilon s}r^{-2\alpha}\mathrm{e}^{-r}\,dr\leq C_{\alpha}\lambda_{0}^{-1},
\]
where $\lambda_{0}$ is the smallest eigenvalue of $A$. Hence $a$ can be chosen
as $C_{\alpha}\lambda_{0}$, for a suitable $C_{\alpha}$.

\noindent\underline{\emph{Step 4}.} We conclude the proof: we have that
$-\mathrm{Tr}[\log(1-2\tilde{a}\tilde{Q}_{s})]\le C_{\alpha}\mathrm{Tr}[\tilde{Q}_{s}]$
since $a$ is small enough, and, as in \emph{step 3},
$$
Tr[Q_{s}]
=   \sum_{\lambda\in\sigma(A)}\lambda^{-2+2\alpha}(2\varepsilon)^{-(1-2\alpha)}\int_{0}^{2\lambda\varepsilon s}r^{-2\alpha}\mathrm{e}^{-r}\,dr
\le C_{\alpha}\varepsilon^{-(1-2\alpha)},
$$
where the sum in $\lambda$ converges since $\alpha<\frac{1}{4}$ and
$\lambda_{n}\approx{n^{\frac{2}{3}}}$. Hence, by \eqref{e:Zchebi} and
\eqref{e:Zserie},
\begin{align*}
\mathbb{P}[\Theta_{\varepsilon}\ge K]
&\le\mathrm{e}^{-\frac{aK^{2}}{\varepsilon}}\mathbb{E}\Bigl[\tilde{a}\bigl[\int_0^1\!\!|AY(s)|^{10}\bigr]^{\frac15}
      +\frac{\tilde{a}^2}2\bigl[\int_0^1\!\!|AY(s)|^8\bigr]^{\frac12}
      +\int_0^1\!\!\mathrm{e}^{\tilde{a}|AY(s)|_H^2}\Bigr]\\
&\le C_\alpha\mathrm{e}^{-\frac{aK^2}\varepsilon}(\mathrm{e}^{C_\alpha\varepsilon^{-(1-2\alpha)}}
      +\varepsilon^{-(1-2\alpha)}
      +\varepsilon^{-2(1-2\alpha)})\\
&\le C_{\#}\mathrm{e}^{-\eta_{\#}\frac{K^2}{\varepsilon}},
\end{align*}
where $\eta_{\#}$ and $C_{\#}$ can be easily found, since $K\geq\frac{1}{2}$.
\end{proof}
\subsection{The deterministic equation}
The basic ingredient of our approach is the bunch of regular paths that every
weak solution has for a positive local (random) time, when the initial
condition is regular. It was called \emph{regular jet} in Flandoli \cite{F}.
It is based on the solutions of the following deterministic equation
\begin{equation}\label{e:determNS}
u(t)+\int_{0}^{t}\left(  Au\left(  s\right)  +B\left(  u,u\right)  \right)
\,ds =x+w\left(  t\right)  .
\end{equation}
We say that
$$
u\in C([0,\infty;H_{\sigma})\cap L_{loc}^2([0,\infty);V)
$$
is a weak solution of \eqref{e:determNS} if
\begin{align*}
\left\langle u(t),\varphi\right\rangle +\int_{0}^{t}\left(  \left\langle
u\left(  s\right)  ,A\varphi\right\rangle -\left\langle B\left(  u\left(
s\right)  ,\varphi\right)  ,u\left(  s\right)  \right\rangle \right)  \,ds
=\left\langle x,\varphi\right\rangle +\left\langle w(t),\varphi\right\rangle
\end{align*}
for every $\varphi\in\mathcal{D}^{\infty}$. Notice that all terms in the above
definition are meaningful, included the quadratic one in $u$ due to the
estimate
$$
\left|  \left\langle B\left(  u,v\right)  ,z\right\rangle \right|  \leq
C\left|  Dv\right|  _{L^{\infty}}\left|  u\right|  _{L^{2}}\left|  z\right|
_{L^{2}}.
$$
We take $w\in\Omega^{\ast}$ where
$$
\Omega^*=\bigcap_{\substack{\beta\in(0,\frac12)\\\alpha\in(0,\frac34)}}C^\beta([0,\infty);D(A^\alpha)).
$$

Consider also the auxiliary Stokes equations
\[
z(t)+\int_{0}^{t}Az\left(  s\right)  \, ds=w\left(  t\right)
\]
having the unique mild solution
\[
z(t)=e^{-tA}w\left(  t\right)  -\int_{0}^{t}Ae^{-(t-s)A}\left(  w\left(
s\right)  -w\left(  t\right)  \right)  \, ds.
\]
From elementary arguments based on the analytic estimates $\left|  A^{\alpha
}e^{-tA}\right|  \leq\frac{C_{\alpha,T}}{t^{\alpha}}$ for $t\in\left(
0,T\right)  $, we have (see for instance Flandoli \cite{Fl1} for details)
$$
z\in C([0,\infty);D(A)).
$$
Let us set
\begin{equation}
\theta_T=\sup_{t\in[0,T]}|Az(t)|.
\end{equation}
Let $C_{0}>0$ be the constant of inequality \eqref{e:inequalityAB} and let
\begin{equation}\label{e:defCstar}
C^{\ast}:=4C_{0}^{2}.
\end{equation}
\begin{lemma}\label{l:primolemmadeterm}
Given $x\in D(A)$ and $w\in\Omega^{\ast}$, let
$K\geq\left|  Ax\right|  $ and $\varepsilon>0$ be such that
$$
(K^2+\theta_\varepsilon^2)(\frac1{2K^{2}}+C^*\varepsilon)<1
$$
Then there exists a solution $u\in C\left(  \left[  0,\varepsilon\right]
;D(A)\right)  $, which is unique in the class of weak solutions, and $\left|
Au\left(  s\right)  \right|  <2K$ for $s\in\left[  0,\varepsilon\right]  $.
%Moreover,
%there is at least one weak solution
%\[
%u\in C\left([0,\infty);H_{\sigma}\right)\cap L_{loc}^{2}\left([0,\infty);V\right).
%\]
%Given $T>0$, if for a weak solution we have
%\[
%\sup_{t\in\left[  0,T\right]  }\left|  Au\left(  t\right)  \right|  <\infty
%\]
%then there exists a unique solution $u\in C\left(  [0,T];D(A)\right)  $.
%Finally, if for a given $T>0$ there is a solution $u\in C\left(
%[0,T];D(A)\right)  $, then it is unique on $[0,T]$ also in the class of weak solutions.
%
%
\end{lemma}
\begin{proof}
We show only the quantitative estimate, the other statements being standard in
the theory of Navier-Stokes equations. For simplicity, all computations will
be made on the limit problem, although they should be made on its Galerkin
approximations. The uniqueness of local solution ensures that the procedure is
nevertheless correct.

Set $v=u-z$, then
\[
\frac{dv}{dt}+Av+B\left(  u,u\right)  =0
\]
and, by using \eqref{e:inequalityAB},
\begin{align*}
\frac{d}{dt}\left|  Av\right|  ^{2}+2\left\|  Av\right\|  _{V}^{2} &
\leq2\left|  \left\langle Av,AB\left(  u,u\right)  \right\rangle \right|
\leq2\left\|  Av\right\|  _{V}\left|  A^{1/2}B\left(  u,u\right)  \right|  \\
&  \leq2C_{0}\left\|  Av\right\|  _{V}\left|  Au\right|  ^{2}\leq\left\|
Av\right\|  _{V}^{2}+C_{0}^{2}\left|  Au\right|  ^{4}\\
&  \leq\left\|  Av\right\|  _{V}^{2}+C^{\ast}(|Av|^{2}+|Az|^{2})^{2}.
\end{align*}
Hence on $[0,\varepsilon]$ we have that
\[
\frac{d}{dt}\left|  Av\right|  ^{2}\leq C^{\ast}(|Av|^{2}+\theta_{\varepsilon
}^{2})^{2},
\]
and so, if we set $y(t)=\left|  Av\left(  t\right)  \right|  ^{2}%
+\theta_{\varepsilon}^{2}$, it follows that
\[
\frac{dy}{dt}\leq C^{\ast}y^{2},\qquad\text{on }\left[  0,\varepsilon\right]
.
\]
Consequently, since $y>0$ (except for the irrelevant case $w\equiv0$), we
have
\[
y(s)\leq\frac{y(0)}{1-C^{\ast}sy(0)},
\]
namely,
\[
\left|  A\left(  u\left(  s\right)  -z\left(  s\right)  \right)  \right|
^{2}+\theta_{\varepsilon}^{2}\leq\frac{\left|  Ax\right|  ^{2}+\theta
_{\varepsilon}^{2}}{1-C^{\ast}s\left(  \left|  Ax\right|  ^{2}+\theta
_{\varepsilon}^{2}\right)  }%
\]
for $s\in\left[  0,\varepsilon\right]  $. Therefore
\[
\left|  Au\left(  s\right)  \right|  ^{2}\leq\frac{2(\left|  Ax\right|
^{2}+\theta_{\varepsilon}^{2})}{1-C^{\ast}s\left(  \left|  Ax\right|
^{2}+\theta_{\varepsilon}^{2}\right)  }\leq\frac{2(K^{2}+\theta_{\varepsilon
}^{2})}{1-C^{\ast}s\left(  K^{2}+\theta_{\varepsilon}^{2}\right)  }.
\]
This result is true until $1-C^{\ast}s\left(  K^{2}+\theta_{\varepsilon}%
^{2}\right)  >0$, namely for $s\in\lbrack0,\frac{1}{C^{\ast}\left(
K^{2}+\theta_{\varepsilon}^{2}\right)  })$. The assumption of the lemma
ensures that $\left[  0,\varepsilon\right]  $ is included in this interval.
Thus the last inequality is true at least on $\left[  0,\varepsilon\right]  $.
Moreover, again by the assumption of the lemma,
\[
\frac{2(K^{2}+\theta_{\varepsilon}^{2})}{4K^{2}}<1-C^{\ast}s\left(
K^{2}+\theta_{\varepsilon}^{2}\right)
\]
that implies
\[
\frac{2(K^{2}+\theta_{\varepsilon}^{2})}{1-C^{\ast}s\left(  K^{2}%
+\theta_{\varepsilon}^{2}\right)  }<4K^{2},
\]
and thus $|Au(s)|^{2}<4K^{2}$, for $s\in\lbrack0,\varepsilon]$.
\end{proof}

\begin{corollary}\label{c:corollaryappendix}
Assume there are $K>0$ and $\varepsilon>0$ such that
$$
\varepsilon\le\frac1{5C^*K^2}
\qquad\text{and}\qquad
\theta_\varepsilon^2\le\frac{1}{4}K^{2},
$$
then, for every $x\in D(A)$ such that $\left|  Ax\right|  \leq K$, we have
$\left|  Au\left(  s\right)  \right|  <2K$ for $s\in\left[  0,\varepsilon
\right]  $.
\end{corollary}
\bibliographystyle{amsplain}

\end{document}